\documentclass{elsart}
\usepackage{amsmath}
\usepackage{epsfig,amsmath,amsfonts}

\newtheorem{theorem}{Theorem}[section]

\usepackage{algorithm,algorithmic}

\begin{document}
\begin{frontmatter}
\title{\Large\bf Efficient computation of high index Sturm-Liouville eigenvalues for problems in physics}
\author{V. Ledoux\thanksref{Foot}\corauthref{cor}},
\author{M. Van Daele and G. Vanden Berghe}
\address{Vakgroep Toegepaste Wiskunde en Informatica, Ghent University, Krijgslaan 281-S9, B-9000 Gent, Belgium}
\thanks[Foot]{Postdoctoral Fellow of the Fund for Scientific Research
- Flanders (Belgium) (F.W.O.-Vlaanderen)}
\corauth[cor]{Corresponding author, email Veerle.Ledoux@UGent.be}
\maketitle
\begin{abstract}
Finding the eigenvalues of a Sturm-Liouville problem can be a computationally
challenging task, especially when a large set of eigenvalues
is computed, or just when particularly large eigenvalues are sought. This
is a consequence of the highly oscillatory behaviour of the solutions corresponding
to high eigenvalues, which forces a naive integrator to take
increasingly smaller steps. We will discuss some techniques that yield
uniform approximation over the whole eigenvalue spectrum and can take
large steps even for high eigenvalues. In particular, we will focus on methods based on coefficient approximation which replace the coefficient functions of the Sturm-Liouville problem by simpler approximations and then solve the approximating problem. The
use of (modified) Magnus or Neumann integrators allows to extend the coefficient approximation idea to higher order methods.
\end{abstract}
\begin{keyword}Sturm-Liouville, Schr\"odinger, eigenvalue, Magnus series, Neumann series
\end{keyword}
\end{frontmatter}
\section{Introduction}
The classical Sturm-Liouville problem (SLP) consists of a linear second-order ordinary differential equation written in formally self-adjoint form
\begin{equation}\label{eq1}
-(p(x)y')'+q(x)y=\lambda w(x)y,
\end{equation}
defined over an interval $a<x<b$ with appropriate boundary conditions at $a$ and $b$. 
An {\em eigenvalue} is a value of $\lambda$ for which \eqref{eq1} has a nontrivial solution $y$ subject to the boundary conditions, and the solution, unique up to scalar multiples, is the associated {\em eigenfunction}. 

The {\em regular} theory which dates back to Sturm (1809-1882) and Liouville (1803-1855) assumes that the coefficient functions are well-behaved, say $p(x), q(x)$ and $w(x)$ are piecewise continuous with $p$ and $w$ strictly positive, on a bounded closed interval $[a,b]$, and that regular boundary conditions are imposed, namely
\begin{equation}\label{bcs}
a_1y(a)+a_2p(a)y'(a)=0,\quad b_1y(a)+b_2p(a)y'(a)=0,
\end{equation}
where $a_1,a_2$ are not both zero, nor are $b_1,b_2$. Then there is an infinite sequence of eigenvalues
\[\lambda_0<\lambda_1<\lambda_2<\dots\]
and eigenfunctions
\[y_0(x),y_1(x),y_2(x),\dots,\]
such that $y_k(x)$ has just $k$ zeros on the open interval $(a,b)$, and such that distinct eigenfunctions are orhogonal with respect to the weight $w(x)$:
\begin{equation}
\int_a^by_i(x)y_j(x)w(x)dx=0,\quad i\neq j.
\end{equation}

The mathematical theory of SLPs is immense (see e.g. \cite{Zettl2005}) but they are not just objects of interest to mathematicians alone. 
Since the early 19th century SLPs have been ubiquitous in applied mathematics since they are the one-dimensional models of a large number of important physical processes in fields such as acoustics, geophysics, waveguide theory, hydrodynamic stability, neutron transport, \dots.  They arise in the analysis of such processes in more than one dimension by the method of Separation of Variables. Another reason why SLPs are of vital interest to physicists is that Schr\"odinger's equation in one dimension is of Sturm-Liouville form.

Some eigenvalue problems have explicit solutions, and are therefore important in the analytical investigation of different physical models. However most eigenvalue problems are not solvable, and computationally efficient approximation techniques are of great applicability.
The numerical solution of (regular) Sturm-Liouville problems is not trivial. The choice of numerical method for efficiently approximating a sequence of eigenvalues of the SLP depends on the desired accuracy of the estimates and also upon the number of eigenvalues required. General ODE boundary-value software can solve SLPs, but inefficiently. 
The challenges are to do this more cheaply, especially when long runs of higher-order eigenvalues are required. \\In fact, many classical methods involve the approximation of the corresponding eigenfunctions by piecewise polynomials and are thus inefficient for the computation of higher eigenvalues which have severely oscillatory eigenfunctions.
There have e.g.\ been many developments in the basic approach of reduction to a matrix eigenproblem using finite differences and finite elements. 
Excellent surveys of such matrix methods are given in \cite{Andrew1994,Andrew2000,VandenBerghe2007}.  Matrix methods can only be used to approximate the first few eigenvalues. The error for even moderately large $k$ is considerable unless the dimension of the associated matrix is very large. 
In this paper we will concentrate on a different class of methods which are based on shooting-type algorithms. These methods perform much better than the matrix methods for singular (or nearly singular) problems, for the computation of eigenfunctions,
and even for the highly accurate computation of the first eigenvalues.
We will discuss in particular some important contributions to the efficient and accurate computation of the higher eigenvalues of SLPs.

\section{Shooting methods}
Shooting methods are based on the reduction of the boundary value problem \eqref{eq1}-\eqref{bcs} to the solution of an initial value problem. The differential equation is solved as an initial value problem over the range $[a,b]$ for a succession of trial values of $\lambda$ which are adjusted till the boundary conditions at both ends can be satisfied at once, at which point we have an eigenvalue. The simplest technique is to `shoot' from $a$ to $b$. This means that one chooses initial conditions which satisfy the boundary condition \eqref{bcs} in $a$:
 \[ y_L(a)=-a_2, \quad p(a)y'_L(a)=a_1\] 
 The boundary condition at $b$ determines `target' values; if the value of $y$ matches the target, we have found an eigenvalue.
In fact, the eigenvalues are determined, using some iterative technique, as the solution of 
\[\phi(\lambda)=b_1y_L(b,\lambda)+b_2p(b)y'_L(b,\lambda)=0.\]
Another option is to shoot from the two ends to some {\em matching point} $a<x_m<b$.
In this case also a right-hand solution $y_R$ is defined satisfying the conditions
\[y_L(b)=-b_2, \quad p(b)y'_L(b)=b_1.\]
The two solutions $y_L$ and $y_R$ are arbitrarily normalised, so their values can always be made to agree at the matching point by renormalising them. However, the criterion for a trial value of $\lambda$ to be an eigenvalue is that the first derivatives should match, as well as te values. The {\em mismatch function} (also called miss-distance) is thus given by
\[\phi(\lambda)=y_L(x_m,\lambda)p(x_m)y'_R(x_m,\lambda)-y_R(x_m,\lambda)p(x_m)y'_L(x_m,\lambda).\]
Numerically some choices of $x_m$ may make it more difficult to compute $\phi(\lambda)$ than do others (see \cite{Pryce1993}). Generally it is a good idea to take the matching point in the interior of the interval, away from singular endpoints.

Thus the procedure for finding the numerical value of an eigenvalue,
consists in evaluating the mismatch function $\phi(\lambda)$, numerically, and then through a finite
series of iterations finding the value of $\lambda$ such that $\phi(\lambda) = 0$ to the required 
accuracy. The usual iterative methods for finding the roots of a function may be
employed here to find the zeros of $\phi(\lambda)$.
One problem associated with this approach is that the function $\phi(\lambda)$ does not give any way of determining
the index of the eigenvalue once it has been found. Thus we have no way of knowing
which eigenvalue we have found when $\phi(\lambda) = 0$. Likewise, in order to converge on a specific eigenfunction, one has to enhance the algorithm, for instance by counting the zeros of the solution as part of the integration for each trial $\lambda$ value.

To solve this problem, the polar coordinate substitution, known as the Pr\"ufer transformation, is used. This transformation makes it possible to specify the eigenvalue to be computed. 
The (scaled) Pr\"ufer transformation is defined by the equations
\begin{equation}\label{polarscaled}y = S^{-1/2}\rho\sin \theta,\quad py' = S^{1/2}\rho\cos \theta,\end{equation}
where $S$ is a strictly positive `scaling function' chosen to give good numerical behaviour. In \cite{Pryce1993} it shown that the resulting differential equations for $\rho$ and $\theta$ are then of the form
\begin{eqnarray}
\theta'&=&\frac{S}{p}\cos^2\theta+\frac{(\lambda w -q)}{S}\sin^2\theta+\frac{S'}{S}\sin\theta\cos\theta,\label{thetaeqp}\\
\frac{2\rho'}{\rho}&=&\left(\frac{S}{p}-\frac{(\lambda w -q)}{S}\right)\sin 2\theta-\frac{S'}{S}\cos 2\theta.\label{eqstr}
\end{eqnarray}
The regular boundary conditions \eqref{bcs}
at $a$ and $b$ define the conditions for $\theta$
\[ \theta(a)=\alpha,\;\;\theta(b)=\beta,\]
where
\[ \tan \alpha =-\frac{S(a)a_2}{a_1},\;\;\tan \beta =-\frac{S(b)b_2}{b_1}.\]
These equations only determine $\alpha$ and $\beta$ up to a multiple of $\pi$. As stated by the following theorem (proved in \cite{Pryce1993}), each (appropriate) choice of this multiple specifies in fact precisely one eigenvalue. 
\begin{theorem}\label{theoremth}Consider the scaled Pr\"ufer equations of a regular Sturm-Liouville problem. 
Let the boundary values $\alpha$ and $\beta$ satisfy the following normalization:
\begin{equation}\label{normalzation}\alpha\in[0,\pi),\quad\beta\in(0,\pi].\end{equation}
Then the $k$th eigenvalue is the value of $\lambda$ giving a solution of \eqref{thetaeqp} satisfying 
\[\theta(a;\lambda)=\alpha,\quad\theta(b;\lambda)=\beta+k\pi.\]
\end{theorem}
This leads to several useful numerical methods based on some form of the Pr\"ufer transformation. Pr\"ufer based shooting methods can be constructed where the counting of the zeros of $y(x)$ needed to compute the specific eigenvalue with a given index $k$ is built in. 
One can for instance define a shooting method for the $\theta$ equation: 
For any $\lambda$, let $\theta_L(x;\lambda)$ and $\theta_R(x;\lambda)$ then be the solutions of \eqref{thetaeqp} satisfying
\begin{equation}\label{condsab}\theta_L(a;\lambda)=\alpha\in[0,\pi),\quad\theta_R(b;\lambda)=\beta\in(0,\pi].\end{equation}
The scaled Pr\"ufer mismatch function is then defined by
\begin{equation}\label{schmiss}\phi(\lambda)=\theta_L(x_m;\lambda)-\theta_R(x_m;\lambda).\end{equation}
and the eigenvalue $\lambda_k$ is the unique value such that $\phi(\lambda_k)=k\pi$.

The SLEIGN code (and its successor SLEIGN2) from Sandia Laboratories \cite{SLEIGN,Bailey2001} uses an (explicit) Runge-Kutta method to integrate
the $\theta$ equation. For certain problems, where a good scaling function $S$ is heuristically found, the oscillations are removed and large steps can be taken. However there is no general method for finding a good scaling function and such shooting methods based on standard initial value libraries often suffer from stepsize restriction when solving for large eigenvalues. They
also have some difficulties caused by stiffness of the $\theta$ equation \eqref{thetaeqp} in a `barrier' region where $(\lambda w-q)/p$ is large and negative.  Instead of using a standard initial value library code, it is a better idea to combine a Pr\"ufer formulation with coefficient approximation, in which the coefficient functions are piecewisely approximated by low degree polynomials. Then the integrations may be performed analytically and stiffness is no longer a problem.

\section{Coefficient approximation}
An important class of methods for the numerical solution of Sturm-Liouville problems is based on coefficient approximation. The basic idea here is to replace the coefficient functions $p(x)$, $q(x)$, $w(x)$ of the Sturm-Liouville equation piecewisely by low degree polynomials so that the resulting equation can be solved analytically.
The idea dates back at least to Gordon \cite{Gordon1969} and Canosa and De Oliveira \cite{Canosa1970} and was studied also by Ixaru \cite{Ixaru1972}, Paine and de Hoog \cite{Paine1980} and Smooke \cite{Smooke1982}. But the standard reference is due to Pruess  \cite{Pruess1973,Pruess1975}. He examined the piecewise constant case and his strategy has been implemented by Pruess and Fulton in the code SLEDGE \cite{SLEDGE}. 
This so-called Pruess method replaces the SLP \eqref{eq1} by the approximating problem
\begin{equation}\label{approxprob}-({\bar p}{y}'(x))'+{\bar q}{y}(x)=\lambda{\bar w}{y}(x),\;\;\;x\in(a,b)\end{equation}
where ${\bar p}, {\bar q}$, and ${\bar w}$ are piecewise constant (midpoint) approximations of the functions $p, q$, and $w$.
The ${y}(x)$ of the approximating problem \eqref{approxprob} can then be integrated explicitly in terms of trigonometric and hyperbolic functions: Let $p$, $q$ and $w$ have constant values ${\bar p}_i$, ${\bar q}_i$, ${\bar w}_i$ in the $i$th interval 
$(x_{i-1}, x_i),\,i=1,\dots,n$ with step size $h_i=x_i-x_{i-1}$:
\[-({{\bar p}_i}{y}'(x))'+{{\bar q}_i}{y}(x)=\lambda{{\bar w}_i}{y}(x)\]
the solution over $[x_{i-1},x_i]$ is then advanced by the relation
\begin{equation}\label{propalgo}
\left(\begin{matrix}
y(x_i)\\{\bar p}_iy'(x_i)
\end{matrix}\right)=
\left(\begin{matrix}
\xi(Z_i)&h_i\eta_0(Z_i)\\Z_i\eta_0(Z_i)/h_i&\xi(Z_i)
\end{matrix}\right)
\left(\begin{matrix}
y(x_{i-1})\\{\bar p}_iy'(x_{i-1})
\end{matrix}\right)
\end{equation}
with $Z_i=h_i^2({\bar q}_i-\lambda {\bar w}_i)/{\bar p}_i$ and
\begin{equation}
\xi(Z) = \left\{ \begin{array}{ll}
\displaystyle \cos(|Z|^{1/2}) & {\rm\ if\ } Z \leq 0\,,\\[.2cm]
\displaystyle \cosh(Z^{1/2}) & {\rm\ if\ } Z > 0\,,
\end{array} \right. \quad
\eta_0(Z) = \left\{ \begin{array}{ll}
\displaystyle \sin(|Z|^{1/2})/|Z|^{1/2} & {\rm\ if\ } Z < 0\,,\\[.2cm]
\displaystyle 1 & {\rm\ if\ } Z=0\,,\\[.2cm]
\displaystyle \sinh(Z^{1/2})/Z^{1/2} & {\rm\ if\ } Z > 0\,,
\end{array} \right.
\label{a2}\end{equation}
One can also propagate the solution from $x_i$ to $x_{i-1}$, by taking the inverse of the transfer matrix in \eqref{propalgo}, which is just the result of replacing $h_i$ by $-h_i$ in this matrix.\\
This gives us a method for explicitly integrating $(y, py')$ over the $x$ range, and to use a
shooting method. This is done e.g. in SLEDGE and combined with the ideas
based on the Pr\"ufer substitution to be able to home in on a particular eigenvalue (see \cite{SLEDGE}).

Pruess proved that if ${\bar p}$, ${\bar q}$ and ${\bar w}$ are piecewise constant on a mesh of typical meshsize $h$ and equal to $p$, $q$, and $w$ at the mesh midpoints and if ${\bar \lambda}_k$ is the $k$th eigenvalue of the approximating problem, then
\[|\lambda_k-{\bar \lambda}_k|\leq Ch^2k|\lambda_k|\]
for all $k$ and small enough $h$.
Thus one would expect a higher eigenvalue to need more meshpoints to compute to a given relative tolerance than a lower eigenvalue. However, as mentioned in \cite{Pryce1993,Marletta1992}, there are two reasons why this is not seen in practice. Firstly, arguments show also that 
$|\lambda_k-{\bar \lambda}_k|\leq Ch|\lambda_k|$ for large $k$ and small enough $h$.
 Secondly, many problems occur in Liouville normal form (Schr\"odinger form) where $p=w=1$ and for these there is an improved error bound \[{|\lambda_k-{\bar \lambda}_k|}\leq Ch^2k^{-1}|\lambda_k|.\]
Thus we can actually use larger $h$ for large $k$ for a given relative error. 

The Pruess-type methods have some important advantages. As already noted, Pruess methods are relatively unaffected by the stiffness/instability which can force a very small stepsize on a standard initial-value solver, and the accuracy is maintained (or even improved) as $k\to\infty$.  Further advantages of the Pruess methods are that it allows a very simple interval truncation algorithm for singular problems (see section \ref{secsing}) and that unlike a method based on a standard initial-value solver, it is practical
with the Pruess method to fix the mesh and evaluate the coefficient midpoint values once
for all before the start of the shooting process. Since the overall shooting process consists of a number of integrations with different
values of $\lambda$ the latter can give a big speed advantage.

A drawback of the Pruess methods is the difficulty in obtaining higher order methods.
It is usual to implement them using Richardson extrapolation. It is clear that the step sizes must be sufficiently small such that the error introduced by the approximation by piecewise constants is not too large. This means that for problems with strongly varying coefficient functions the number of intervals in a mesh can be quite large. Some approaches have been suggested towards the realization of higher order methods based on coefficient approximation. These approaches can be classified as ``modified integral series methods'', which will be discussed next. 


\section{Modified integral series methods}
We will consider two integral series which allow the natural extension of the Pruess-ideas to higher order methods: a Neumann series and a Magnus series. In fact, these integral series offer an easy way to approximate the coefficient functions of the SLP by higher order (piecewise) polynomials, giving more accurate results than the approximation by a piecewise constant.

\subsection{The Neumann and Magnus expansion}
There is an emerging family of numerical methods based on integral series representation of ODE solutions. Consider the linear differential equation
\begin{equation} {\bf y}'=A(x){\bf y},\;\;{\bf y}(0)={\bf y}_0\in {\mathbb R}^N.\label{Mfirsteq}\end{equation}
The simplest integral series is obtained by applying Picard iteration \cite{Hairer1987} to obtain the fundamental solution of the matrix linear ODE
\begin{equation} \begin{split}{\bf y}(x)=&\Big[I+\int_0^xA(x_1)d{x_1}+\int_0^xA(x_1)\int_0^{x_1} A(x_2)dx_2d{x_1}
\\&+\int_0^xA(x_1)\int_0^{x_1} A(x_2)\int_0^{x_2}A(x_3)dx_3dx_2d{x_1}+\dots\Big]{\bf y}_0\end{split}\end{equation}
This series is known as the Feynman-Dyson path ordered exponential in quantum mechanics, in mathematics it is known as the {\em Neumann} series or Peano series.\index{Neumann series}\index{Peano series}

The Magnus and Cayley expansions are two other examples. They are obtained by transforming Eq.\ \eqref{Mfirsteq} to the suitable Lie algebra\index{Lie algebra} and applying the Picard iteration to the transformed ODE. Details on both approaches can be found in \cite{Iserles2000}. The Cayley expansion is based on the Cayley transform while the Magnus expansion is based on the exponential map. The approach of Magnus \cite{Magnus1954} aims at writing the solution of Eq.\ \eqref{Mfirsteq} as
\[ {\bf y}(x)=\exp(\Omega(x)){\bf y}_0\] where
$\Omega(x)$ is a suitable matrix. The {\em Magnus expansion} says that
\begin{equation}\label{omegas}\begin{split}\Omega(x)=&\int_0^xA(x_1)d{x_1}-\frac{1}{2}\int_0^x\left[\int_0^{x_1} A(x_2)dx_2,A(x_1)\right]d{x_1}
\\&+\frac{1}{4}\int_0^x\left[\int_0^{x_1}\left[\int_0^{x_2}A(x_3)dx_3,A(x_2)\right]dx_2,A(x_1)\right]d{x_1}
\\&+\frac{1}{12}\int_0^x\left[\int_0^{x_1} A(x_2)dx_2,\left[\int_0^{x_1} A(x_3)dx_3,A(x_1)\right]\right]d{x_1}+\dots
\end{split}\end{equation}
where $[\cdot,\cdot]$ denotes the matrix commutator\index{commutator} defined by $[X,Y]=XY-YX$.

Numerical schemes based on the Magnus expansion received a lot of attention due to their preservation of Lie group symmetries (see \cite{Iserles2000,Iserles1999} and references therein). The Neumann series does not respect Lie group structure but avoids the use of the matrix exponential. The use of Neumann series integrators has been proved successfull for certain large, highly oscillatory systems in 
\cite{Iserles2004}.

Since the SLP can be written in the matrix form \eqref{Mfirsteq}, both Neumann and Magnus schemes can be considered for the numerical solution of the SLP. The Sturm-Liouville equation \eqref{eq1} in matrix form reads
\begin{equation}\label{matrixSLP}{\bf y}'(x)=A(x){\bf y}(x)=\left(\begin{matrix}0&1/p(x) \\q(x)-\lambda w(x)&0\end{matrix}\right) {\bf y}(x)\end{equation}
with ${\bf y}^T=(y(x),p(x)y'(x))$. Moan \cite{Moan1998} was the first to consider a Magnus series integrator for the SLP in the Schr\"odinger form $y''(x)=(q(x)-\lambda)y(x)$ or in matrix form
\begin{equation}\label{matrixschrod}{\bf y}'(x)=\left(\begin{matrix}0&1 \\q(x)-\lambda &0\end{matrix}\right) {\bf y}(x).\end{equation}
He applied the Magnus integrator directly to this problem. However poor approximations were obtained for the higher eigenvalues, as a result of
the finite radius of convergence of the Magnus series \cite{Niesen}. 
When the solution of a linear system ${\bf y}'=A(x){\bf y}$ oscillates rapidly, modified schemes should be used, as recommended in \cite{Iserles2002,Iserles2002b,Degani2006}. 
Describing these modified schemes we will focus on the basic Schr\"odinger equation $y''(x)=(q(x)-\lambda)y(x)$, but the schemes can be extended to the more general Sturm-Liouville problem $-(p(x)y'(x))'+q(x)y(x)=\lambda w(x)y(x)$.

\subsection{Modified Neumann and Magnus schemes for the Schr\"odinger equation}
We consider the Sturm-Liouville problem in Schr\"odinger form eq.\ \eqref{matrixschrod} which is a problem of the form
\begin{equation}\label{orig}
{\bf y}(x)'=A(x,\lambda){\bf y}(x),\;\;{\bf y}(a)={\bf y}_0,
\end{equation}
where ${\bf y}=[y(x), y'(x)]^{T}$. Note that the coefficient matrix is in $sl(2)$,i.e.\ the matrix has a zero trace.

Suppose that we have already computed ${\bf y}_{i-1}\approx{\bf y}(x_{i-1})$ and that we wish to advance the numerical solution to $x_{i}=x_{i-1}+h_i$. We first compute a constant approximation ${\bar q}$ of the potential function $q(x)$
\begin{equation}\label{constapproxq}
{\bar q}=\frac{1}{h_i}\int_{x_{i-1}}^{x_{i-1}+h_i}q(x)dx.\end{equation}
Next we change the frame of reference by letting
\begin{equation}\label{transfo}
{\bf y}(x)=e^{(x-x_{i-1}){\bar A}}{\bf u}(x-x_{i-1}),\quad x_{i-1}\leq x\leq x_{i}
\end{equation}
where
\begin{equation}
{\bar A}(\lambda)=\left(\begin{matrix}0&1\\{\bar q}-\lambda&0\end{matrix}\right).
\end{equation}
We treat ${\bf u}$ as our new unknown which itself obeys the linear differential equation
\begin{equation}\label{modeq}
{\bf u}'(\delta)=B(\delta,\lambda){\bf u}(\delta),\quad \delta\in[0,h_i],\quad {\bf u}(0)={\bf y}_{i-1}
\end{equation}
where
\begin{equation}
B(\delta,\lambda)=e^{-\delta{\bar A}}\left(A(x_{i-1}+\delta)-{\bar A}\right)e^{\delta{\bar A}}.
\end{equation}
The matrix $B$ can be computed explicitly. With $\xi(Z)$ and $\eta_0(Z)$ defined as in eq.\ \eqref{a2}
we can write $B$ as
\begin{equation}\label{oscB}
B(\delta,\lambda)=\Delta_q(\delta)\left(\begin{matrix}\delta\eta_0(Z_{2\delta})
&\displaystyle\frac{1-\xi(Z_{2\delta})}{2(\lambda-{\bar q})}
\\\displaystyle-\frac{1+\xi(Z_{2\delta})}{2}&
-\delta\eta_0(Z_{2\delta})\end{matrix}\right),
\end{equation}
where $\Delta_q(\delta)={\bar q}-q(x_{i-1}+\delta)$ and $Z_{\gamma}=Z(\gamma)=({\bar q}-\lambda)\gamma^2$.

We have thus replaced one linear system by another. The new system \eqref{modeq} has one crucial advantage over \eqref{orig}: the entries of the matrix $B$ are themselves rapidly oscillating functions (for $\lambda >{\bar q}$). This is not very helpful when \eqref{modeq} is solved by a classical method as e.g.\ a Runge-Kutta or multistep method. When the modified equation is however solved by an integral series method, repeated evaluation of integrals of $B$ is required. This integration is a ``smoothing'' operator: the amplitude is decreased once
the integrand is integrated. As a result the higher the oscillation, the faster the convergence of the integral series method and the faster the decay in local error. We can refer to \cite{Iserles2002,Iserles2002b,Degani2006} for numerical results confirming the success of this approach for highly oscillatory ODEs.

Over each interval $[x_{i-1} , x_{i}]$ an integral series is applied on the transformed equation ${\bf u}'(\delta) = B(\delta){\bf u}(\delta)$. This requires the truncation of the integral series and the replacement of integrals by quadrature (see next section).
The solution ${\bf y}$ in $x = x_{i}$ is then obtained from ${\bf y}(\lambda, x_{i}) = e^{h_i {\bar A}}{\bf u}(h_i )$. 
Note that $e^{h_i{\bar A}}$ is the known solution of the system with constant potential
\begin{equation}
{\rm expm}\left(\begin{matrix}0&h_i\\h_i({\bar q}-\lambda) &0\end{matrix}\right)=\left(\begin{matrix}\xi(Z_h)&h_i{\eta_0(Z_h)}\\
Z_h{\eta_0 (Z_h)}/h_i&\xi(Z_h)\end{matrix}\right),\quad Z_h=Z(h_i)
\end{equation}
and thus the same as the transfer matrix in \eqref{propalgo}. 

The first option we consider is the use of a Neumann scheme.
Application of the Neumann series integrator to the modified equation ${\bf u}'(\delta) = B(\delta){\bf u}(\delta)$ gives
\begin{equation}\label{Neushrod}\begin{split}{\bf u}(h_i)=\displaystyle \Big[I+\int_0^{h_i}B(x){\rm d}x+\int_0^{h_i}\int_0^{x_1} B(x_1)B(x_2){\rm d}x_2{\rm d}x_1
\\\displaystyle \hspace*{10mm}+\int_0^{h_i}\int_0^{x_1} \int_0^{x_2}B(x_1)B(x_2)B(x_3){\rm d}x_3{\rm d}x_2{\rm d}x_1+\dots\Big]{\bf y}_{i-1}\end{split}\end{equation}
When only the first term in the Neumann series is retained, one has ${\bf u}(h_i)={\bf y}_{i-1}$ and with ${\bf y}(x_{i}) = e^{h_i {\bar A}}{\bf u}(h_i )$
this is exactly the second-order Pruess method given by eq.\ \eqref{propalgo}. Higher order methods are obtained by including more Neumann terms. In \cite{Degani2006} it was shown that in fact each extra Neumann term can be seen as a correction term in a Piecewise constant Perturbation Method (PPM) of Ixaru and co-workers.  The PPMs use a perturbation technique (the first systematic description of this technique is due to Ixaru \cite{Ixaru1984}) to construct some correction terms which are added to the known solution of the approximating problem $y''=({\bar q}-\lambda)y$ with a piecewise constant potential ${\bar q}$. The more correction terms included the higher the order of the algorithm. 
In \cite{Ixaru1997,Ixaru2000} the PPM algorithm is described and applied on regular Schr\"odinger problems. The PPMs were also extended to general Sturm-Liouville problems using the Liouville transform and formed the basis of the Fortran code SLCPM12 \cite{SLCPM12} and the graphical Matlab software package  {\sc Matslise} \cite{Ledoux2005a}. \\
To approximate the integrals in \eqref{Neushrod} quadrature must be used which can deal adequately with the oscillatory entries of the matrix function $B$. In section \ref{Filonsection} a Filon-type quadrature rule will be discussed which is very similar to the procedure used in the description of high order PPMs in \cite{Ixaru1997,ledoux2004}. There the potential function $q$ is replaced by a piecewise polynomial, which makes the integrals in \eqref{Neushrod} analytically solvable. The degree of the (piecewise) polynomial can be taken sufficiently large such that this approximation has no influence on the accuracy of the method. The order of the method then only depends on the number of terms retained in the Neumann expansion (i.e.\ PPM correction terms). As mentioned, including only the first Neumann term gives us a method of order two. Including also the second Neumann term (the first integral) leads to a method of order four, adding the third term (the double integral) results in an eighth order method and adding the fourth term gives us a method of order ten (see \cite{Iserles2004}). 

Another option is to apply a Magnus method to the modified equation \eqref{modeq}. The Magnus expansion is then 
\begin{equation}
\label{expansion}\sigma(\delta)=\sigma_1(\delta)+\sigma_2(\delta)+\sigma_3(\delta)+\sigma_4(\delta)+\dots,\end{equation}
where
\begin{eqnarray}\label{Magnschrod}
\nonumber \sigma_1(\delta)&=&\int_0^\delta B(x)dx,\\
\nonumber \sigma_2(\delta)&=&-\frac{1}{2}\int_0^\delta\int_0^{x_1}[B(x_2),B(x_1)]dx_2d{x_1},\\
\nonumber \sigma_3(\delta)&=&\frac{1}{12}\int_0^\delta\left[\int_0^{x_1}B(x_2)dx_2,\left[\int_0^{x_1}B(x_2)dx_2,B(x_1)\right]\right]d{x_1},\\
\sigma_4(\delta)&=&\frac{1}{4}\int_0^\delta\left[\int_0^{x_1}\left[\int_0^{x_2}B(x_3)dx_3,B(x_2)\right]dx_2,B(x_1)\right]d{x_1},
\end{eqnarray}
and ${\bf u}(\delta)=e^{\sigma(\delta)}{\bf y}_{i-1}, \,\delta\geq 0$. Thus, to compute ${\bf y}_{i}=e^{h{\bar A}}e^{\sigma(h)}{\bf y}_{i-1}$ with $h=h_i$, we need to approximate $\sigma(h)$ by truncating the expansion \eqref{expansion} and replacing integrals by quadrature (see \ref{Filonsection}). As shown in \cite{Degani2006}, truncating all but the first integral leads to a fourth order method, while including also $\sigma_2$ gives us a scheme of order eight. Having approximated $\sigma(h)$, its $2 \times 2$ matrix exponential must be computed.
We note that ${\sigma(h)}$ is always a two by two matrix with zero trace. For such matrices the following is true:
\begin{equation}\label{expmzerotrace}
{\rm expm}\left(\begin{matrix}a&b\\c &-a\end{matrix}\right)=\left(\begin{matrix}\xi(\omega)+a\eta_0(\omega)&b\eta_0(\omega)\\c\eta_0(\omega)&\xi(\omega)-a\eta_0(\omega)\end{matrix}\right),\quad \omega={a^2+bc}.
\end{equation}
Here $a,b,c,\omega$ are functions of $x$ and $E$. 


\subsection{Quadrature of the (multivariate) integrals}
\label{Filonsection}
Practical implementation of both the Neumann and Magnus series requires the replacement of multivariate integrals by quadrature. Although multivariate quadrature is usually considered a hard problem, it is possible to implement Neumann and Magnus expansions with surprisingly cheap and effective quadrature. Moreover when a Filon-type quadrature method is used, even the highly oscillating integrals, which appear when $\lambda \gg {\bar q}$, are approximated to a suitable precision in a small number of function evaluations per step. Filon quadrature has been analysed extensively in \cite{Iserles2004b}.

For a Neumann integrator as well as for a Magnus integrator, the univariate (modified) integral $\int_0^{h_i} B(\delta)d\delta$ needs to be approximated. A Filon-type rule is used. Here this means that $\Delta_q(\delta)$ in \eqref{oscB} is replaced by a polynomial, i.e.\ by the Lagrange polynomial \begin{equation}\label{Lagr}{\mathfrak L}_{\Delta_q}(\delta)=\sum_{l=1}^\nu \Delta_q(c_l h_i)\ell_l(\delta)\end{equation} where $\ell_l$ is the $l$th cardinal polynomial of Lagrangian interpolation and $c_1$, $c_2$, $\dots$,$c_\nu$ are distinct quadrature nodes. The resulting integrals can then be solved analytically. For each entry in the univariate integral a scheme of the following form results 
\[ h_i\sum_{l=1}^{\nu}b_l(\omega)\Delta_q(c_lh_i),\quad \omega={\bar q}-\lambda.\]
When no further Neumann or Magnus terms are retained in the algorithm, the truncated Neumann or Magnus scheme is of order four and it is then sufficient to have $\nu=2$ Legendre quadrature nodes (or $\nu=3$ Lobatto nodes). This means thus that in this case $\Delta_q$ is approximated by a linear polynomial. \\
For schemes of order eight, the double integral must be included and $\nu=4$ Legendre nodes should be used. As for the univariate integral, the double integral is computed by replacing $\Delta_q $ by the polynomial ${\mathfrak L}_{\Delta_q}$ and solving the resulting integrals analytically (using a symbolic software package). Each entry in the double integral is then approximated by an expression of the form
\[h_i^2\sum_{k=1}^{\nu}\sum_{l=1}^{\nu}b_{k,l}(\omega)\Delta_q(c_kh_i)\Delta_q(c_lh_i),\]
where the values of $\Delta_q$ that have been already evaluated for
the quadrature of the univariate integral are reused.
For triple and further integrals the same procedure can be applied: replace $\Delta_q$ by the Lagrange polynomial of sufficiently high degree and then use the resulting analytic expressions for the integrals as approximating formulae.
Note that also the value of ${\bar q}$ in \eqref{constapproxq} is computed by Gauss-Legendre with $\nu$ nodes, and thus the same function evaluations of $q(x)$ are needed as to compute the different $\Delta_q(c_lh_i)$ in the Lagrange polynomial.

An alternative way to apply the Filon-type rule is by approximating $q(x)$ (piecewisely) by a series over shifted Legendre polynomials (as in done in the description of PPM \cite{Ixaru1997,ledoux2004}):
\begin{equation}
q(x)\approx\sum_{s=0}^{\nu-1}Q_sh_i^sP_s^*(\delta/h_i),\quad \delta=x-x_{i-1}
\end{equation}
By the method of least squares the expressions for the coefficients $Q_s$ are obtained:
\begin{eqnarray}
{Q}_s&=&\frac{(2s+1)}{h_i^{s+1}}\int_0^h q(x_{i-1}+\delta)P^*_s(\delta/h_i) d\delta,\;\;m=0,1,2,\dots. 
\label{Gauss}
\end{eqnarray}
It can then be noted that ${\bar q}=Q_0$ and $\Delta_q(\delta)\approx{\mathfrak L}_{\Delta_q}(\delta)=-\sum_{s=1}^{\nu-1}Q_sh_i^sP_s^*(\delta/h_i)$.
Writing ${\mathfrak L}_{\Delta_q}(\delta)$ in this form is fully equivalent as using \eqref{Lagr}, but allows to obtain shorter expressions for the formulae approximating the integrals.
To compute the integrals \eqref{Gauss} Gauss-Legendre is used, requiring $\nu$ function evaluations of $q$. Suppose we truncate all but the first integral, resulting in a method of order four. We need to discretise the integral consistently with the order of the method. To this end, we take $\nu=2$.
With 
 \[{\widehat{\xi}}=\xi(Z_{2h}),\quad {\widehat \eta_0}=\eta_0(Z_{2h}),\quad Z_{2h}=4Z_h=4({\bar q}-\lambda)h_i^2\]
and ${\hat Q}_s=h_i^{s+1}Q_s,s=1,\dots,\nu-1$,
we then obtain the following
\begin{eqnarray}\label{UniFilon}
\nonumber \frac{1}{h_i}\int_0^{h_i} \Delta_q(\delta)\delta\eta_0(Z_{2\delta})d\delta&\approx& \frac{{\hat Q}_1(-1-{\widehat{\xi}}+2{\widehat \eta_0})}{4Z_h}\\
 \nonumber\int_0^{h_i} \Delta_q(\delta)\left(1+\xi(Z_{2\delta})\right)d\delta&\approx&\int_0^{h_i} \Delta_q(\delta)\xi(Z_{2\delta})d\delta\\
\nonumber\int_0^{h_i} \Delta_q(\delta)\left(1-\xi(Z_{2\delta})\right)d\delta&\approx&-\int_0^{h_i} \Delta_q(\delta)\xi(Z_{2\delta})d\delta\approx{\hat Q}_1{\widehat \eta_0}+\frac{{\hat Q}_1(1-{\widehat{\xi}})}{2Z_h}
\end{eqnarray}
which allows us to approximate $\int_0^{h_i}B(\delta)d\delta$.  Note that the quadrature approximation of the non-oscillating integral $\int_0^{h_i}\Delta_q(\delta)d\delta$ vanishes: $\int_0^{h_i}{\mathfrak L}_{\Delta q}(\delta)d\delta=\int_0^{h_i} ({\bar q}-\mathfrak{L}_{q}(\delta)) d\delta=0$, since both $\int_0^{h_i} \mathfrak{L}_{q}(\delta) d\delta$ and ${\bar q}h_i$ are the Gauss-Legendre quadrature approximations of $\int_0^{h_i} q(x_{i-1}+\delta)d\delta$.

To construct a method of order eight, we need to take $\nu=4$ and we have to include the double integral. To compute the univariate integral we have now
\begin{eqnarray}
\nonumber \frac{1}{h_i}\int_0^{h_i} \Delta_q(\delta)\delta\eta_0(Z_{2\delta})d\delta&\approx& \frac{({\hat Q}_1+3{\hat Q}_2+6{\hat Q}_3){\widehat \eta_0}}{2Z_h}-\frac{({\hat Q}_3+{\hat Q}_1)({\widehat{\xi}}+1)+{\hat Q}_2({\widehat{\xi}}-1)}{4Z_h}\\\nonumber&&+\frac{3{\hat Q}_2(1-{\widehat{\xi}})-15{\hat Q}_3({\widehat{\xi}}+1)}{4Z_h^2}+\frac{15{\hat Q}_3{\widehat \eta_0}}{2Z_h^2}\\
\nonumber-\int_0^{h_i} \Delta_q(\delta)\xi(Z_{2\delta})d\delta&\approx&({\hat Q}_1+{\hat Q}_2+{\hat Q}_3){\widehat \eta_0}+\frac{(3{\hat Q}_2+15{\hat Q}_3){\widehat \eta_0}}{Z_h}\\\nonumber&&-\frac{3{\hat Q}_2({\widehat{\xi}}+1)+({\hat Q}_1+6{\hat Q}_3)({\widehat{\xi}}-1)}{2Z_h}+\frac{15{\hat Q}_3(1-{\widehat{\xi}})}{2Z_h^2}.
\end{eqnarray}Suppose we construct a Magnus method, we consider then the approximation of $\sigma_2$. A similar procedure can be followed to compute the double integral in a Neumann method.
As in \cite{Ledoux2008} we write the double integral in $\sigma_2$ as
\begin{equation}\label{integr3}
\begin{split}
\int_0^{h_i}\int^{\delta_1}_0[B(\delta_2),B(\delta_1)]d\delta_2d\delta_1=2\int_0^{h_i}\int^{\delta_1}_0 \Delta_q(\delta_1)\Delta_q(\delta_2)K_1(\delta_1,\delta_2)d\delta_2d\delta_1 U_1\\+
2\int_0^{h_i}\int^{\delta_1}_0 \Delta_q(\delta_1)\Delta_q(\delta_2)K_2(\delta_1,\delta_2)d\delta_2d\delta_1 U_2\\+2\int_0^{h_i}\int^{\delta_1}_0 \Delta_q(\delta_1)\Delta_q(\delta_2)K_3(\delta_1,\delta_2)d\delta_2d\delta_1 U_3
\end{split}
\end{equation}
where $K_1(x,y)=\displaystyle{y}\eta_0(Z_{2y})-x\eta_0(Z_{2x})$, $K_2(x,y)=\xi(Z_{2x})-\xi(Z_{2y}),$ 
$K_3(x,y)=\displaystyle{(x-y)}\eta_0(Z_{2(x-y)})$ and
\begin{equation}\label{US}\nonumber
U_1=\left(\begin{matrix}0&\frac{1}{2(\lambda-{\bar q})}\\\frac{1}{2}&0\end{matrix}\right),\quad\quad U_2=\left(\begin{matrix}-\frac{1}{4(\lambda-{\bar q})}&0\\0&\frac{1}{4(\lambda-{\bar q})}\end{matrix}\right),\quad\quad U_3=\left(\begin{matrix}0&\frac{1}{2(\lambda-{\bar q})}\\\frac{-1}{2}&0\end{matrix}\right).
\end{equation}
The three integrals in \eqref{integr3} are replaced by quadrature by again approximating $\Delta_q$ by the polynomial ${\mathfrak L}_{\Delta_q}$ .
The expression for the third integral is then for instance:
\begin{equation}\label{sigma2}\nonumber
\begin{split}
&\int_0^{h_i}\int_0^{\delta_1}\Delta_q(\delta_1)\Delta_q(\delta_2)K_3(\delta_1,\delta_2)d\delta_2d\delta_1\approx\\
&\hspace*{5mm}\Big(\frac{{\hat Q}_2^2-{\hat Q}_3^2-{\hat Q}_1^2-2{\hat Q}_3{\hat Q}_1}{4Z_h}+\frac{-{\hat Q}_1^2+15{\hat Q}_2^2-66{\hat Q}_3^2-42{\hat Q}_3{\hat Q}_1}{4Z_h^2}\\&\hspace*{5mm}+\frac{9{\hat Q}_2^2-405{\hat Q}_3^2-30{\hat Q}_3{\hat Q}_1}{4Z_h^3}-
\frac{225{\hat Q}_3^2}{4Z_h^4}\Big){\widehat \eta_0}+\Big(\frac{{\hat Q}_1^2-3{\hat Q}_2^2+6{\hat Q}_3^2+7{\hat Q}_3{\hat Q}_1}{4Z_h^2}\\&\hspace*{5mm}+
\frac{30{\hat Q}_3{\hat Q}_1-9{\hat Q}_2^2+105{\hat Q}_3^2}{4Z_h^3}+
\frac{225{\hat Q}_3^2}{4Z_h^4}\Big){\widehat{\xi}}-\frac{42{\hat Q}_2^2+70{\hat Q}_1^2+30{\hat Q}_3^2}{840Z_h}-\frac{5{\hat Q}_3{\hat Q}_1}{4Z_h^2}.\end{split}
\end{equation}
In practice, one should use a truncated series expansion for small $Z_h$ values (see \cite{Ledoux2008}).

Even higher order algorithms can be constructed including more Magnus (or Neumann) terms in the scheme. In \cite{Ledoux2008b} a Magnus scheme of order 10 is described where $\nu=5$ and in \cite{ledoux2004} PPM-schemes up to order 16 are presented. Note that only the terms where the degree in $h$ is smaller or equal to the required degree of the method have to be included in the algorithm, for instance in the approximation of $\sigma_3$ and $\sigma_4$ the term in ${\hat Q}_4^3$ can be disregarded in the Magnus scheme of order 10.


The modified Magnus methods and modified Neumann methods are well suited for the repeated solution of the initial value problems which appear in the shooting procedure. These initial value problems are solved for a fixed potential $q$ but for different values of $\lambda$. As shown in \cite{Ixaru1997,Ledoux2008}, an $\lambda$-independent mesh can be computed which is then (re)used in all eigenvalue computations. Moreover also the value ${\bar q}$ and the coefficients $Q_s$ are computed and stored once for all before the start of the shooting process. 
Algorithm \ref{shootalgo} shows the basic shooting procedure in which a modified Magnus algorithm
is used to propagate the left-hand and right-hand solutions.
\begin{algorithm}
\caption{A Sturm-Liouville solver based on a modified Magnus method}\label{shootalgo}
\begin{algorithmic}[1]
\STATE Use stepsize selection algorithm to construct mesh $a=x_0<x_1<...<x_n=b$
\FOR{$i=1$ to $n$}
\STATE Compute ${\bar q}$ and $Q_s,s=1,\dots,\nu-1$ for the $i$th interval (Gauss-Legendre with $\nu$ nodes).
\ENDFOR
\STATE Choose a meshpoint $x_{m}$ ($0\leq m\leq n $) as the matching point.
\STATE Set up initial values for ${\bf y}_L$ satisfying the BC at $a$ and initial values for ${\bf y}_R$ satisfying the BC at $b$. Choose a trial value for $\lambda$.
\REPEAT 
\FOR{$i=0$ to $m-1$} 
\STATE ${\bf y}_L(x_{i+1})=e^{h_i{\bar A}}e^{\sigma(h_i)}{\bf y}_L(x_i)$ 
\vspace*{1mm}
\ENDFOR
\FOR{$i=n$ down to $m+1$}
\STATE ${\bf y}_R(x_{i-1})=e^{-\sigma(h_i)}e^{-h_i{\bar A}}{\bf y}_R(x_i)$ 
\vspace*{1mm}
\ENDFOR
\STATE Adjust $\lambda$ by comparing ${\bf y}_L(x_m)$ with ${\bf y}_R(x_m)$.
\UNTIL{$\lambda$ sufficiently accurate}
\end{algorithmic}
\end{algorithm}




\section{Some notes on singular problems}
\label{secsing}
When the problem is singular, either because $(a,b)$ is an infinite interval or because at least one of the coefficients $p^{-1}$, $q$, $w$ is not integrable up to one of the endpoints, then an interval truncation procedure must be adopted. 
Different algorithms are implemented in the available SLP library codes to determine a truncated endpoint and appropriate boundary conditions to give a prescribed accuracy (see \cite{Pryce1993}). The SLEDGE package \cite{SLEDGE} even has algorithms for automatically classifying the nature of the problem, regular or singular, limit-circle singularity or limit-point singularity and so on. This classification information is important to determine whether or not there is a continuous spectrum, when there are eigenvalues and how many, and what boundary condition should be imposed at a singular endpoint.

As mentioned before, the algorithm applied in the SLEDGE package relies on Pruess coefficient approximation by piecewise constants (namely the midpoint values of each interval), and uses repeated extrapolation to achieve accuracy. In a first pass a crude initial mesh is chosen by an equidistribution process. SLEDGE then repeatedly bisects this initial mesh and uses iterated extrapolation. 
An infinite endpoint is transformed to zero by the (local) change of variable $t=1/x$, and subsequent bisections near the endpoint are done in terms of the variable $t$.
SLEDGE's approach automatically regularizes singular endpoints: evaluating the coefficients at the mesh midpoints can be regarded as truncating the interval at the midpoints of the first and final intervals of the mesh. Every time the mesh is bisected these implicit truncation points move closer to the singular endpoints. The boundary conditions are always applied in the original endpoints.

For the higher order coefficient approximation methods, a similar approach as in SLEDGE can be used to truncate a singular problem. For these methods the coefficients are only evaluated in the Legendre nodes. Since the first and last Legendre node in an interval are not equal to the beginpoint or endpoint of that interval, a singular problem is implicitly truncated. By decreasing the size of the first interval (if $a$ is a singular endpoint) or the last interval (if $b$ is a singular endpoint) in the mesh, the implicit truncation points move closer to the singular endpoint. However no software based on higher order modified integral methods is currently availabe which can automatically handle singular endpoints. For this reason, SLEDGE is still the package of choice for many applications. 

\section{Some experiments}
We consider two well-known test problems. The Coffey-Evans equation is one of the more difficult test problems in the literature (test problem 7 in \cite{Pryce1993}, introduced in \cite{CoffeyEvans}).
It is a Schr\"odinger equation with $q (x) =-2\beta\cos(2x)+\beta^2\sin^2(2x)$  and $y(-\pi/2)=y(\pi/2)=0$ as boundary conditions. The first 50 eigenvalues for $\beta=30$ have been determined. For this potential $\lambda_0$ is close to zero and there are very close eigenvalue triplets $\{\lambda_2,\lambda_3,\lambda_4\}$, $\{\lambda_6,\lambda_7,\lambda_8\},\ldots$ as $\beta$ increases.  
The second problem is a problem from chemical physics: the Woods-Saxon problem \cite{VandenBerghe1989} defined by
\begin{equation}\label{potWS}
q(x)=-50\frac{1-\frac{5t}{3(1+t)}}{1+t}
\end{equation}
with $t=e^{(x-7)/0.6}$ over the interval $[0,15]$. The eigenvalue spectrum of this Woods-Saxon problem contains 14 eigenenergies $\lambda_{0},...,\lambda_{13}$.

Tables \ref{tab:CoffeyEvansProblem} -- \ref{tab:WSProblem2} show some results for the two test problems obtained with five different coefficient approximation methods. The first method used is the Pruess method of order $P=2$. This method is compared with some higher order methods which were discussed in section \ref{Filonsection}: a Neumann method and a Magnus method of order $P=4$, and a Neumann method and  Magnus method of order $P=8$. We present for each problem, a selection of the considered exact eigenvalues $\lambda_k$, and the (absolute) error for the corresponding eigenvalues calculated with the coefficient approximation methods. For the moment equidistant meshes are used in order to allow easier comparison between the different algorithms. An automatic stepsize selection algorithm will be discussed afterwards.


\begin{small}
\begin{table}
	\caption{Absolute value of (absolute) errors $\Delta \lambda_k$ in the eigenvalues computed for the Coffey-Evans problem with different coefficient approximation methods and $n=128$ steps in the equidistant mesh. $a$E-$b$ means $a.10^{-b}$. $P$ is the order of the method.}
	\label{tab:CoffeyEvansProblem}
	\centering
		\begin{tabular}{ll|c|cc|cc}
		\hline
		&&$P=2$&\multicolumn{2}{c}{$P=4$}\vline&\multicolumn{2}{c}{$P=8$}\\
		$k$&$\lambda_k$&&Neumann&Magnus&Neumann&Magnus\\
		\hline
	\phantom{1}0&\phantom{306}0.0000000000000000&1.7E-1&1.3E-3&1.3E-3&6.4E-9&1.0E-7\\
			\phantom{1}1&\phantom{3}117.9463076620687587&1.5E-1&3.5E-3&3.5E-3&1.5E-8& 2.5E-7\\
			\phantom{1}2&\phantom{3}231.6649292371271088&1.3E-1&3.1E-3&3.1E-3&1.6E-9&7.9E-8\\
			\phantom{1}3&\phantom{3}231.6649293129610125&1.3E-1&3.1E-3&3.1E-3&3.4E-8&1.6E-7\\
			\phantom{1}4&\phantom{3}231.6649293887949167&1.3E-1&3.1E-3&3.1E-3&2.3E-9&2.3E-7\\
			\phantom{1}5&\phantom{3}340.8882998096130157&1.0E-1&6.3E-3&6.3E-3&1.5E-8&2.5E-7\\
			\phantom{1}6&\phantom{3}445.2830895824354620&7.7E-2&5.6E-3&5.6E-3&1.0E-8&1.8E-7\\
			\phantom{1}8&\phantom{3}445.2832550313310036&7.7E-2&5.4E-3&5.4E-3&1.0E-8&1.8E-7\\
			10&\phantom{3}637.6822498740469991&3.1E-2&6.7E-3&6.7E-3&5.0E-10&5.7E-9\\
			15&\phantom{3}802.4787986926240517&2.2E-2&5.1E-3&5.1E-3&6.3E-9&9.1E-8\\
			20&\phantom{3}951.8788067965913828&4.6E-2&4.2E-3&4.2E-3&5.6E-9&8.1E-8\\
			30&1438.2952446408023577\phantom{+}&2.3E-2&3.7E-3&3.7E-3&2.5E-9&2.8E-8\\
			40&2146.4053605398535082\phantom{+}&1.3E-2&3.0E-3&3.0E-3&1.1E-9&1.5E-8\\
			50&3060.9234915114205911\phantom{+}&8.9E-3&2.2E-3&2.2E-3&4.5E-10&8.8E-9\\
				\hline
		\end{tabular}
\end{table}
\end{small}

\begin{small}
\begin{table}
	\caption{Absolute value of (absolute) errors $\Delta \lambda_k$ in the eigenvalues computed for the Woods-Saxon problem with different coefficient approximation methods and $n=64$ steps in the equidistant mesh. $a$E-$b$ means $a.10^{-b}$. $P$ is the order of the method.}
	\label{tab:WSProblem}
	\centering
		\begin{tabular}{ll|c|cc|cc}
		\hline
		&&$P=2$&\multicolumn{2}{c}{$P=4$}\vline&\multicolumn{2}{c}{$P=8$}\\
		$k$&$\lambda_k$&&Neumann&Magnus&Neumann&Magnus\\
		\hline
	\phantom{1}0&-49.45778872808258&1.7E-3&5.5E-6&5.5E-6&2.8E-10&4.9E-9\\
			\phantom{1}1&-48.14843042000639&5.1E-3&4.0E-5&4.0E-5&2.3E-9&4.0E-8\\
			\phantom{1}2&-46.29075395446623&9.1E-3&1.3E-4&1.4E-4&9.5E-9&1.6E-7\\
			\phantom{1}3&-43.96831843181467&1.3E-2&3.2E-4&3.2E-4&2.7E-8&4.7E-7\\
			\phantom{1}4&-41.23260777218090&1.8E-2&6.0E-4&6.0E-4&6.1E-8&1.1E-6\\
			\phantom{1}5&-38.12278509672854&2.1E-2&1.0E-3&1.0E-3&1.1E-7&2.0E-6\\
			\phantom{1}6&-34.67231320569997&2.5E-2&1.5E-3&1.5E-3&1.8E-7&3.2E-6\\
			\phantom{1}7&-30.91224748790910&2.7E-2&2.1E-3&2.1E-3&2.6E-7&4.4E-6\\
			\phantom{1}8&-26.87344891605993&2.7E-2&2.8E-3&2.8E-3&3.1E-7&5.5E-6\\
			\phantom{1}9&-22.58860225769320&2.6E-2&3.4E-3&3.4E-3&3.2E-7&5.9E-6\\
			10&-18.09468828212811&2.3E-2&4.0E-3&4.0E-3&2.6E-7&5.1E-6\\
			11&-13.43686904026007&1.7E-2&4.4E-3&4.4E-3&1.1E-7&3.1E-6\\
			12&\phantom{1}-8.67608167074520&7.3E-3&4.6E-3&4.6E-3&1.1E-7&1.1E-7\\
			13&\phantom{1}-3.90823248120989&5.9E-3&4.3E-3&4.3E-3&3.2E-7&3.5E-6\\
				\hline
		\end{tabular}
\end{table}
\end{small}

The five different methods were first applied on the same mesh, with $128$ steps for the Coffey-Evans problem and $64$ steps for the Woods-Saxon problem. Results are shown in tables \ref{tab:CoffeyEvansProblem} and \ref{tab:WSProblem}. All methods allow the approximation of higher eigenvalues or large batches of eigenvalues. However it is clear that the higher order methods need much less mesh intervals to reach a prescribed accuracy. This can also been seen from the tables \ref{tab:CEProblem} and \ref{tab:WSProblem2}, where the number of intervals in the equidistant mesh ($nint$) and function evaluations ($nfev$) are shown that each method needs to reach an accuracy of (approximately) $10^{-8}$. 
The data reported in the four tables enable several conclusions:
\begin{itemize}
\item Each of the five methods can reproduce accurate results, even for high eigenvalues.
\item However when Neumann or Magnus terms are introduced, less intervals are needed in the mesh to reach a certain input accuracy.
\item As a result the number of function evaluations is also decreasing with increasing order of the method. 
The second order Pruess method only needs one function evaluation per interval, while the fourth order methods need two and the eighth order methods need four. However to reach the same accuracy, the second order method needs many more intervals which causes a very high total number of function evaluations. 
\item Since only a small number of function evaluations is needed, the construction of the mesh takes less time for a higher order method. But also the shooting process in which the equation is repeatedly integrated at various values of $\lambda$ is considerably faster as a result of the smaller number of intervals needed in the mesh. The timings shown in tables \ref{tab:CEProblem} and \ref{tab:WSProblem2} were obtained using Matlab-implementations of the different methods.
\item The modified Neumann and Magnus methods are particularly well suited to compute (large) batches of eigenvalues. A remarkably small number of mesh intervals is sufficient to be able to compute accurate approximations for the higher eigenvalues. There is also a big speed advantage in the fact that the repeatedly asked task of integrating the equation at various trial values for an eigenvalue is completely separated from the time-consuming process of constructing a mesh where all function evaluations are performed and saved for later use.
\item The fourth order Magnus and Neumann versions reach the same accuracy. The eighth order Neumann method seems to be somewhat more precise than its Magnus counterpart. This is a consequence of the finite radius of convergence of the Magnus expansion (see \cite{Niesen}) which means that the steps in the ${\bar q}\gg\lambda$ region can not always be taken as large as for the Neumann expansion when $Z_h$ is large and positive. Note that cases with very large positive $Z$ rarely appear in practice. They are e.g.\ ruled out by WKB arguments which usually shrink the interval \cite{Ledoux2006inf}. 
\item The Neumann methods are somewhat faster than their Magnus counterpart. The evaluation of the matrix exponential in the Magnus method requires little extra time for the second order Sturm-Liouville problem. For problems with higher dimension, the computation of the matrix exponential may be fairly expensive.
\item However truncated Neumann expansions do not respect Lie-group structure which can be a limitation in some applications. But as remarked in \cite{Iserles2004}, the basic step underlying the Neumann expansion discussed here is the transformation \eqref{transfo}, which always preserves Lie-group structure. Departures from a Lie group might occur only in the function ${\bf u}$, in other words in the correction term. This results in far less severe loss of Lie-group structure than is the case with classical Runge-Kutta or multistep methods.
\end{itemize}

\begin{small}
\begin{table}
	\caption{Absolute value of (absolute) errors $\Delta \lambda_k$ in the eigenvalues computed for the Coffey-Evans problem. Different coefficient approximation methods were used on an equidistant mesh with $nint$ steps to reach an accuracy of approximately $10^{-8}$. $nfev$ is the number of evaluations of the potential function and $T(s)$ is the CPU-time needed to compute the first 51 eigenvalues.}
	\label{tab:CEProblem}
	\centering
		\begin{tabular}{ll|c|cc|cc}
		\hline
		&&$P=2$&\multicolumn{2}{c}{$P=4$}\vline&\multicolumn{2}{c}{$P=8$}\\
		$k$&$\lambda_k$&&Neumann&Magnus&Neumann&Magnus\\
		\hline
	    \phantom{1}0&\phantom{100}0.0000000000000000&6.7E-7&1.9E-8&1.9E-8&6.3E-8&1.0E-7\\
			10&\phantom{1}637.6822498740469991&1.3E-7&1.1E-7&1.1E-7&8.1E-9&5.6E-9\\
			20&\phantom{1}951.878806796591382&1.7E-7&7.3E-8&7.3E-8&5.7E-8&8.1E-8\\
			30&1438.2952446408023577&8.4E-8&7.4E-8&7.4E-8&2.3E-8&2.8E-8\\
			40&2146.4053605398535082&4.7E-8&7.4E-8&7.4E-8&1.4E-8&1.5E-8\\
			50&3060.9234915114205911&4.5E-8&5.6E-8&5.6E-8&8.4E-9&8.5E-9\\
				\hline
				$nint$&&65536&2048&2048&96&128\\
				$nfev$&&65536&4096&4096&384&512\\
				$T(s)$&&34.14&2.53&4.18&0.85&0.99\\
				\hline
		\end{tabular}
\end{table}
\end{small}

\begin{small}
\begin{table}
	\caption{Absolute value of (absolute) errors $\Delta \lambda_k$ in the eigenvalues computed for the Woods-Saxon problem. Different coefficient approximation methods were used on an equidistant mesh with $nint$ steps to reach an accuracy of approximately $10^{-8}$. $nfev$ is the number of evaluations of the potential function and $T(s)$ is the CPU-time needed to compute the 14 eigenvalues.}
	\label{tab:WSProblem2}
	\centering
		\begin{tabular}{ll|c|cc|cc}
		\hline
		&&$P=2$&\multicolumn{2}{c}{$P=4$}\vline&\multicolumn{2}{c}{$P=8$}\\
		$k$&$\lambda_k$&&Neumann&Magnus&Neumann&Magnus\\
		\hline
	\phantom{1}0&-49.45778872808258&3.3E-9&8.5E-11&8.5E-11&1.1E-11&1.8E-10\\
			\phantom{1}2&-46.29075395446623&1.7E-8& 2.1E-9&2.1E-9&3.8E-10&6.2E-9\\
			\phantom{1}4&-41.23260777218090&3.4E-8&9.6E-9&9.6E-9&2.4E-9&4.0E-8\\
			\phantom{1}6&-34.67231320569997&4.7E-8&2.5E-8&2.5E-8&7.4E-9&1.2E-7\\
			\phantom{1}8&-26.87344891605993&5.3E-8&4.7E-8&4.7E-8&1.3E-8&2.2E-7\\
			10&-18.09468828212811&4.5E-8&7.3E-8&7.3E-8&1.2E-8&2.0E-7\\
			12&\phantom{1}-8.67608167074520&1.6E-8&9.0E-8&9.0E-8&2.3E-9&7.4E-9\\
				\hline
				$nint$&&32768&1024&1024&96&96\\
				$nfev$&&32768&2048&2048&384&384\\
				$T(s)$&&4.27&0.38&0.62&0.19&0.22\\
				\hline
		\end{tabular}
\end{table}
\end{small}

Of course using a uniform mesh is rarely a good idea e.g.\ when dealing with (truncated) singular problems. For automatic software a stepsize selection algorithm should be used. We refer to \cite{SLEDGE} for the procedure used in SLEDGE. A meshing algorithm has also been proposed for the PPM in e.g.\ \cite{Ixaru1997} and for the numerical solution of regular Schr\"odinger problems with an eighth order Magnus method in \cite{Ledoux2008}.
\begin{figure}
	\begin{center}
	\includegraphics[width=11cm]{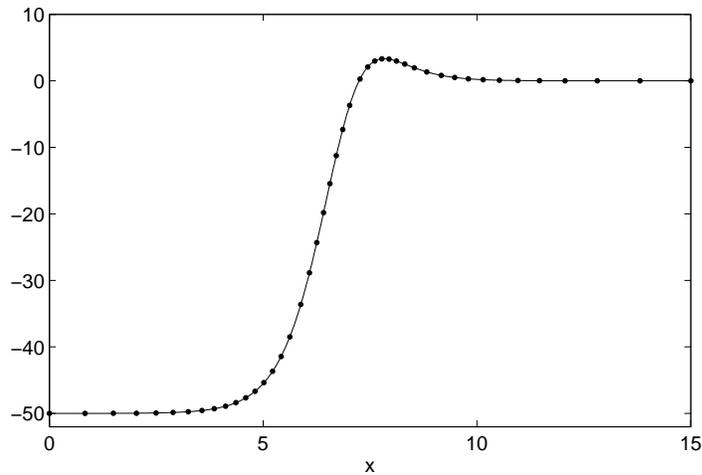}
	\end{center}
	\caption{The mesh resulting from the stepsize selection algorithm of the eighth order Neumann method for the Woods-Saxon potential with input tolerance $tol=10^{-6}$: the values of the potential $q(x)$ at the mesh points are marked by dots.}\label{figpart2}
\end{figure}   
\begin{small}
\begin{table}
	\caption{Absolute value of (absolute) errors $\Delta \lambda_k$ in the eigenvalues computed for the Woods-Saxon problem. An eighth order Neumann method is used on the mesh (with 46 steps) shown in figure \ref{figpart2}.}
	\label{tab:WSProblem4}
	\centering
		\begin{tabular}{ll|c}
		\hline
		$k$&$\lambda_k$&$\Delta \lambda_k$\\
		\hline
	\phantom{1}0&-49.45778872808258&1.5E-10\\
			\phantom{1}2&-46.29075395446623&4.6E-9\\
			\phantom{1}4&-41.23260777218090&5.3E-8\\
			\phantom{1}6&-34.67231320569997&2.4E-7\\
			\phantom{1}8&-26.87344891605993&2.0E-7\\
			10&-18.09468828212811&1.5E-7\\
			12&\phantom{1}-8.67608167074520&1.4E-7\\
				\hline
		\end{tabular}
\end{table}
\end{small}
We briefly describe a procedure which can be used for the eighth order Neumann method and apply it on the Woods-Saxon problem.
The resulting mesh is shown in Figure \ref{figpart2}, and the results obtained over this mesh are listed in Table \ref{tab:WSProblem4}.\\
The stepsize selection algorithm applied here is based on a local error estimate. Let $tol$ be an input tolerance parameter. In order to attain the local error estimate $\epsilon_{i}=tol$ over the interval $i$, the stepsize $h_{i}$ is chosen as a function of the previous stepsize $h_{i-1}$ as follows. First we compute
\begin{equation}{\bar h}_{i}= h\left(\frac{tol}{\epsilon}\right)^{1/8},\quad h=h_{i-1}\end{equation} 
where $\epsilon$ is the error estimate.
A decision is taken in terms of $\Delta=|{\bar h}_{i}/h-1|$. If $\Delta>0.1$ the procedure is repeated with $h={\bar h}_{i}$. If $\Delta \leq 0.1$, $h$ is accepted as the stepsize for the interval $i$.\\
To construct the error estimate $\epsilon$ we consider the maximal (absolute) difference between the entries of the transfer matrix of an (embedded) sixth order scheme and the entries of the transfer matrix of our eighth order algorithm.  This difference is evaluated through the sum of all terms in ${\bar Q}_3$. When scanning over $\lambda$ we have used only three values, those such that $Z({\bar h}_i)=({\bar q}-\lambda){\bar h}_i^2=-m^2\pi^2, m=0,1,2$. The selection of only these was mainly intended to speed up the evaluation but this is enough as the error decreases like $O(1/\sqrt{\lambda})$ for large $Z$-values \cite{Degani2006,Ixaru1997}, and as confirmed by experimental tests showing that the error is indeed larger for smaller values of $Z$. 

Let us consider now a singular problem, just to illustrate that the modified integral series methods can be extended to singular problems. The general Woods-Saxon problem is a Schr\"odinger problem of the form
\[y''(x)=\left(\frac{l(l+1)}{x^2}+q(x)-\lambda\right)y(x)\]
with $q(x)$ as in \eqref{potWS} and $x\in[0,+\infty]$. When the orbital quantum number $l$
equals zero, the potential is a well behaved, nonsingular function, as the Woods-Saxon problem we considered before. When $l>0$ the problem is singular in the origin. We will compute the eigenvalues for the problem with $l=2$. This problem has 13 eigenvalues in its spectrum. 
The infinite endpoint can be dealt with by a change of variable converting (implicitly) to a finite interval (as in SLEDGE) or with explicit interval truncation using WKB arguments as in \cite{Ledoux2006inf}. We describe here only a way to deal with the singularity in the origin and consider the (truncated) problem over the interval [0,20]. We applied the stepsize selection algorithm discussed above over the interval $[\epsilon,20]$ and added to the resulting mesh the interval $[0,\epsilon]$. Note that the potential function is only evaluated in the Legendre points of the interval $[0,\epsilon]$ and not in the singular endpoint $a=0$. A first approximation of an eigenvalue is computed over the mesh. Then the interval $[0,\epsilon]$ is bisected and a further eigenvalue approximation is computed. The process of bisection in $[0,\epsilon]$ is repeated until two successive eigenvalue approximations agree within the user specified tolerance. At each iteration, the shooting algorithm for the next eigenvalue approximation is started using the approximation last obtained. Table \ref{tab:WSProblem5} shows the results for two different values of $\epsilon$. A larger value of $\epsilon$ requires of course more bisections ($nbisec$) of the interval $[0,\epsilon]$. The results obtained are within the requested accuracy $tol=10^{-7}$. 
\begin{small}
\begin{table}
	\caption{Absolute value of (absolute) errors $\Delta \lambda_k$ in the eigenvalues computed for the singular Woods-Saxon problem. An eighth order Neumann method is used with user input tolerance $tol=10^{-7}$. $nint$ is the number of intervals in the initial mesh and $T$ is the CPU-time needed to compute all eigenvalues.}
	\label{tab:WSProblem5}
	\centering
		\begin{tabular}{ll|cc}
		\hline
		$k$&$\lambda_k$&$\epsilon=0.01$&$\epsilon=0.1$\\
		\hline
	\phantom{1}0&-48.349481052120&6.7E-11&8.8E-9\\
			\phantom{1}2&-44.121537377319&6.4E-10&1.5E-8\\
			\phantom{1}4&-38.253426539679&2.1E-9&1.2E-8\\
			\phantom{1}6& -31.026820921773&1.5E-9&2.2E-9\\
			\phantom{1}8&-22.689041510178&1.1E-8&2.4E-8\\
			10&-13.52230335295&5.6E-11&5.8E-8\\
			12&\phantom{1}-3.972491432846&2.3E-8&1.5E-9\\
				\hline
				$nint$&&144&111\\
				$T(s)$&&0.91&1.18\\
				$nbisec$&&1&5\\
				\hline
		\end{tabular}
\end{table}
\end{small}




\section{Conclusion}
In this paper we discussed some techniques which allow the efficient approximation of high eigenvalues of a Sturm-Liouville problem. We focused in particular on algorithms based on approximation of the coefficient functions of the differential equation. The simplest coefficient approximation method is the so-called Pruess method, which replaces the coefficient functions over each mesh interval by the midpoint value in that interval and then solves (analytically) the approximating problem. This Pruess method has some significant advantages over shooting methods based on standard initial-value solvers especially when looking for higher eigenvalues and was implemented in the well-known Sturm-Liouville solver SLEDGE. The oscillations in the eigenfunctions no longer determine (restrict) the step sizes. Now the step sizes depend on the errors made in replacing the coefficient functions by  piecewise constant approximations. It is clear that larger steps could be taken when the coefficient functions are replaced by higher order polynomials. However when the coefficient functions are replaced by polynomials of degree greater than zero, the approximating problem is not really easier to solve than the original problem. Therefore only piecewise constant (and linear) polynomial approximations were used for a long time.
But, Neumann or Magnus integral series offer a way to construct methods based on higher order (piecewise) polynomial approximation. Using (modified) Neumann or Magnus schemes, we can construct methods which still allow the easy analytic integration of an approximating problem with piecewise constant coefficients but use higher order polynomial approximations to construct some extra (correction) terms. Depending on the number of terms included in the Neumann or Magnus series, algorithms of different orders can be constructed. Experiments show that indeed these Neumann and Magnus integrators share the advantages of the Pruess method and allow to approximate high eigenvalues in a remarkably small number of steps. Also a singular problem was considered to illustrate that singular endpoints can be dealt with in a relatively simple way.


\providecommand{\Gr}{Gr} \providecommand{\ea}{et al}

\end{document}